\documentclass{llncs}

\usepackage[final,                     
color,notref,notcite]{showkeys}         

\usepackage{makeidx}  
\usepackage{epsfig}
\begin{document}
\title{SM stability for time-dependent problems}
\titlerunning{SM stability for time-dependent problems}
\author{Petr N. Vabishchevich}
\authorrunning{Petr N. Vabishchevich}
\tocauthor{Petr N. Vabishchevich
(Keldysh Institute of Applied Mathematics, RAS)}
\institute{Keldysh Institute of Applied Mathematics, RAS\\
4 Miusskaya Square, 125047 Moscow, Russia\\
\email{vabishchevich@gmail.com}}
\maketitle

\begin{abstract}
Various classes of stable finite difference schemes 
can be constructed to obtain a numerical solution. It
is important to select among all stable schemes such a 
scheme that is optimal in terms of certain additional
criteria.
In this study, we use a simple boundary value problem 
for a one-dimensional parabolic equation to discuss the 
selection of an approximation with respect to time. 
We consider the pure diffusion equation, the pure convective transport 
equation and combined convection-diffusion phenomena.
Requirements for the unconditionally stable
finite difference schemes are formulated that are related 
to retaining the main features of the differential
problem. The concept of SM stable finite difference scheme 
is introduced. The starting point are difference schemes 
constructed on the basis of the various Pad$\acute{e}$ approximations. 
\end{abstract}

\section{Introduction}

When time-dependent problems of mathematical physics 
are solved numerically, much emphasis is
placed on computational algorithms of higher 
orders of accuracy (e.g., see \cite{Hundsdorfer,Gustafsson}). 
Along with improving
the approximation accuracy with respect to space, 
improving the approximation accuracy with respect to
time is also of interest. In this respect, the 
results concerning the numerical methods for ordinary 
differential equations (ODEs) \cite{Ascher,LeVeque} 
provide an example. 
Taking into account the specific features of time-dependent 
problems for PDEs, we are interested in numerical 
methods for solving the Cauchy problem
in the case of stiff equations \cite{Rakit,Hairer,Butcher}.

When time-dependent problems are solved approximately, 
the accuracy can be improved in various
ways. In the case of two-level schemes (the solution 
at two adjacent time levels is involved), polynomial
approximations of the scheme operators on the solutions 
are used explicitly or implicitly.
The most popular representatives of such schemes 
are Runge-Kutta methods \cite{Butcher,Verwer}, which are widely
used in modern computations. The main feature of 
the multilevel schemes (multistep methods) manifests
itself in the approximation of time derivatives with 
a higher accuracy on a multipoint stencil. A characteristic 
example is provided by multistep methods based on 
backward numerical differentiation \cite{Gear}.

Various classes of stable finite difference schemes 
can be constructed to obtain a numerical solution
\cite{SamTheor,SamVabMat}. It
is important to select among all  stable schemes such
a scheme that is optimal in terms of certain additional
criteria. In the theory of finite difference schemes, 
there is the class of asymptotically stable schemes (see
\cite{SamGul,SamVab}) that ensure the correct 
long-time behavior of the approximate solution. 
In the theory of numerical
methods for ODEs (see \cite{Butcher,Gear}), the concept 
of $L$-stability is used, which reflects the 
long-time asymptotic behavior of the approximate 
solution from a different point of view.

In   \cite{VabSM} there are considered properties of two-level difference schemes of high order approximation for the approximate solution
of the Cauchy problem for evolutionary equations with self-adjoint operators.  
The simplest boundary value problem for the one-dimensional parabolic equation serves as a basic problem.  
The concept of SM stability (Spectral Mimetic stability) of a difference scheme is introduced. 
This property is connected with the behavior of individual harmonics of the approximate solutions.   

In this paper, we continue to study the SM properties of difference schemes for the approximate solution of 
unsteady problems of mathematical physics.  
For the model boundary value problem for one-dimensional parabolic equation there are considered the spectral characteristics of 
approximations in space and in time. 
In particular, good approximation properties are observed for the convection operator of third order in space.  
Features two-level schemes of higher order of approximation in time, based on the Pad$\acute{e}$ approximation, are considered for solving problems of 
mathematical physics with symmetric and skew-symmetric operators.  

\section{Problem formulation}

We consider finite-dimensional real Hilbert space $H$, where the scalar product and the norm is $(\cdot,\cdot)$ and $\|\cdot\|$, respectively.  
Let $u(t)$ ($0 \leq t \leq T > 0$) is defined as the solution of the Cauchy problem for evolutionary equation of first order:  
\begin{equation}\label{1}
  \frac{d u}{d t} + \Lambda \, u = f(t),
  \quad 0 < t \leq T ,
\end{equation}
\begin{equation}\label{2}
  u(0) = u_0 . 
\end{equation}
The right part $f(t) \in H$ of equation  (\ref{1}) is given and $\Lambda$ is a linear non-negative, in general, nonself-adjoint  
operator from $H$ in $H$ ($\Lambda  = \Lambda(t) \geq 0$) depending on  $t$.

For problem (\ref{1}), (\ref{2}) the estimate of stability is easily established. 
Taking into account the skew-symmetric property of operator $\Lambda$, we have the equality   
\[
  \|u\| \frac{d \|u\|}{dt} =   (f, u).
\]
With regard
\[
  (f, u) \le  \|u\| \|f\|
\]
we obtain a simple estimate of stability for the solution of (\ref{1}), (\ref{2}) with respect to the initial data and right hand side:  
\begin{equation}\label{3}
  \|u(t)\| \leq 
  \|u_0\| + \int_{0}^{t} \|f(\theta)\| d \theta .
\end{equation}
We would like to preserve these properties of the differential problem after transition to a discrete analogue of problem (\ref{1}), (\ref{2}).

The main attention in our discussion is given to unsteady boundary value problem  for partial differential equations.  
In this context, we can associate the Cauchy problem  (\ref{1}), (\ref{2}) with application of the method of lines (approximation in space).  
Given the importance for applications, we will conduct our consideration on an example of a boundary value problem 
for the one-dimensional parabolic equation of second order.  
Let a sufficiently smooth function  $u(x,t)$ satisfies the equation  
\begin{equation}\label{4}
  \frac{\partial  u}{\partial  t} +   
  \mathcal{L} u = 0,
  \quad 0 < x < 1, \quad 0 < t \leq T 
\end{equation}
and the initial conditions  
\begin{equation}\label{5}
  u(x,0) = u_0(x),
  \quad 0 < x < 1 .
\end{equation}
The periodicity of the spatial variable is assumed:
\begin{equation}\label{6}
  u(x+1,t) = u(x,t), 
  \quad 0 < t \leq T .
\end{equation}

We associate operator $ \mathcal{L}$ with the convection-diffusion equation, defining  
\begin{equation}\label{7}
  \mathcal{L} = \chi \mathcal{C} + (1-\chi) \mathcal{D}  
\end{equation}
for some constant $0 \leq \chi \leq 1$. 
Here the operators of convective and diffusive transport are defined according  
\begin{equation}\label{8}
  \mathcal{C} u = \frac{\partial  u}{\partial  x},
  \quad \mathcal{D} u = - \frac{\partial^2  u}{\partial  x^2} .
\end{equation}
If $\chi = 1$ equation (\ref{4}) is the convection transport equation whereas if $\chi = 0$ it is the diffusion equation.  

The discrete problem should inherit the main properties of the differential problem.  
In model problem (\ref{4})--(\ref{6}) the skew-symmetric property of operator $\mathcal{C}$
as well as the self-adjoint and non-negative properties of  operator  $\mathcal{D}$ should be preserved
\begin{equation}\label{9}
  \mathcal{C} = - \mathcal{C}^*,
  \quad \mathcal{D} = \mathcal{D}^* \geq 0
\end{equation}
in space $L_2(0,1)$ for functions satisfying (\ref{6}).
Stability of the solution of the corresponding problem (\ref{1}), (\ref{2}) (estimate (\ref{3})) is provided by 
similar properties of the grid analogs of convective and diffusive transport operators.  
In our research, we are concentrating on the spectral characteristic of the solution
(Spectral Mimetic Properties for grid approximations), when considering the behavior of individual harmonics of the approximate solution.  

\section{SM properties of the approximation in space}

Let us introduce a uniform grid with step  $h$:
\[
  x_i = i h, \ i = 0, \pm 1, \pm 2, ...,  \ M h = 1,
\]
\[
  \omega = \{ x_i \ | \ i = 0, 1,  ...,  \ M -1 \} .
\]
We use standard index-free notations of the theory of difference schemes \cite{SamTheor}.
Let $ w = w_i = w (x_i) $, and  for the left, right and central difference derivatives we set  
\[
  \partial_- w = \frac{w_i - w_{i-1}}{h},
  \quad   \partial_+ w  = \frac{w_{i+1} - w_{i}}{h},
\]
\[
  \partial_0 w = \frac{1}{2} (\partial_- w + \partial_+ w ) =
  \frac{w_{i+1} - w_{i-1}}{2h} 
\]
respectively.

After approximation in space we put for problem (\ref{4})--(\ref{6}) the corresponding discrete problem  
\begin{equation}\label{10}
  \frac{d  y}{d  t} + \Lambda  y = 0,    
  \quad x \in \omega ,
  \quad 0 < t \leq T ,
\end{equation}
\begin{equation}\label{11}
  y(x,0) = u_0(x),
  \quad x \in \omega .
\end{equation}
Some key  possibilities in the choice of grid operator  $\Lambda$ with the focus on the properties of the differential 
operator $\mathcal{L}$ should be noted.  

We define Hilbert space $H=L_2(\omega)$ of periodic   grid functions ($y(x+1)=y(x)$) with the inner product and the norm  
\[
  (y,w) = \sum_{x \in \omega} y(x)\,w(x)\,h,
  \quad \|y\|^2 = (y,y) .
\]
To guarantee  stability of solution (\ref{10}), (\ref{11}), operator $\Lambda$  must be non-negative
($\Lambda \geq 0$) in $H$.
Conservatism (neutral stability) communicates directly with the  skew-symmetric property of operator 
$\Lambda$ ($\Lambda = - \Lambda^*$) in  $H$.

For the convection equation  ($\chi = 1$, $\mathcal{L} = \mathcal{C}$)
with  $f = 0$ the norm of the solution of problem (\ref{4})--(\ref{6}) does not change in time:  
\begin{equation}\label{12}
  \|u(t)\| =   \|u_0\| .
\end{equation}
Equality (\ref{12}) reflects the conservation property of the solution (conservation law), the neutral stability of the solution.  

\textit{1. Upwind (directional) approximations of first order.}
To approximate the convective terms (see, e.g.,  \cite{Hundsdorfer,SamVabConv,morton}),
the upwind  first-order approximations are traditionally widely used.
In this case, the grid convection operator  $C$  has the form  
\begin{equation}\label{13}
  C = \partial_- .
\end{equation}
Operator  $C$ defined in  (\ref{13}) is non-negative ($C \geq 0$). 
In this case, for the solution of problem (\ref{10}), (\ref{11})  the estimate 
\begin{equation}\label{14}
  \| y\| \leq \|u_0\| , 
  \quad 0 < t \leq T 
\end{equation}
is true.  

\textit{2. Central-difference approximation.}
Another well-known variant  is to use approximations of second order where
\begin{equation}\label{15}
  C = \partial_0 .
\end{equation}
In this case we have $C = - C^*$ and for problem  (\ref{10}), (\ref{11}), (\ref{15})
holds  
\begin{equation}\label{16}
  \| y\| = \|u_0\| , 
  \quad 0 < t \leq T .
\end{equation}

\textit{3.Upwind second-order approximations.}
When choosing approximations of higher order (second and above) for the convective terms, 
we are trying at least partially to preserve the properties of the first order approximations, 
which are connected primarily with the  monotonicity  (fulfillment of the  maximum principle). 
The most interesting attempts in the class of linear approximations are associated with the use of approximations 
with the upwind differences of second order \cite{Gustafsson,Hirsch}. 
For our problem (\ref{4})--(\ref{6}), we have  
\[
  C y = \frac{3 y_i - 4 y_{i-1} + y_{i-2}}{2 h} .
\]
Using previously introduced operator notations we obtain 
\begin{equation}\label{17}
  C = \partial_- + \frac{h}{2} \partial_- \partial_- .
\end{equation}
Operator $C \geq 0$ and so  for the solution of problem  (\ref{10}), (\ref{11}) estimate (\ref{14}) holds again.

\textit{4. Approximations of third order.}
In computing  practice third order approximations are not in common use.  
In fact, they are only mentioned (see, e.g., \cite{Gustafsson,LeVeque})  without any meaningful analysis.  
In this case the difference convection operator can be written in the form  
\begin{equation}\label{18}
  C = \partial_0 - \frac{h^2}{6} \partial_- \partial_- \partial_+.
\end{equation}
In index notation equation (\ref{18}) takes the form  
\[
  С y = \frac{2 y_{i+1} + 3 y_i - 6 y_{i-1} + y_{i-2}}{6 h} .
\]
Operator $C \geq 0$ and its energy (equal to $(C y,y)$) is three times less energy operator  $C$, 
which is defined  by rule (\ref{17}) (upwind second-order approximations).  

The stability conditions (neutral stability) of the considered approximations of convection transfer are 
associated with the general properties of the operator (non-negativity, skew-symmetric property).  
More detailed information gives us the spectrum of the difference operator, 
its proximity to the spectrum of differential operator. 
This inheriting the properties  of the differential problem for the difference problem at the spectral level 
we associate \cite{VabSM} with SM properties.    

Consider the corresponding differential problem for eigenvalues and eigenfunctions.  
For operator $\mathcal{C}$ we have 
\begin{equation}\label{19}
  \frac{d  v}{d  x} = \lambda \, v,
  \quad 0 < x < 1,
\end{equation}  
\begin{equation}\label{20}
  v(x+1) = v(x) . 
\end{equation}  
The solution of spectral problem (\ref{19}),(\ref{20}) is  
\[
  \lambda_m = i 2 \pi m,
\]
\[
  v_m(x) = e^{i 2 \pi m x},
  \quad m = 0, \pm 1, \pm 2, ... .
\]
For solution of problem (\ref{4})--(\ref{6})  we obtain the representation  
\begin{equation}\label{21}
  u(x,t) = \sum_{m=-\infty}^{\infty} (u_0, v_m) e^{-\lambda_m t} v_m(x),
\end{equation}
where  
\[
  (u_0, v_m) = \int_{0}^{1} u_0(x)v_m(x) \, d x ,
  \quad m = 0, \pm 1, \pm 2, ...\quad \]
are the coefficients of  expansion for function $u_0 (x)$.

We now consider the appropriate spectral problems of the grid problem  
\begin{equation}\label{22}
  C v = \mu v
\end{equation}
with the above-mentioned approximations of convective transport.  
For definiteness, we assume that $M$ is odd.  
Eigenfunctions of problem (\ref{22}) for the considered difference operators (\ref{13}), (\ref{15}), (\ref{17})  and  (\ref{18}) 
have the form  
\begin{equation}\label{23}
  w_m(x) = e^{i 2 \pi m x},
  \quad x \in \omega ,
  \quad m = 0, \pm 1, \pm 2, ..., \frac{M-1}{2}.
\end{equation}

For difference operator  $C$, defined according to  (\ref{13}), the eigenvalues have the form  
\[
  \mu_m = \frac{1 - e^{i 2 \pi m h}}{h},
  \quad m = 0, \pm 1, \pm 2, ..., \frac{M-1}{2}.
\]
For the imaginary and real parts we obtain  
\begin{equation}\label{24}
  \mathrm{Re} \ \mu_m = \frac{2}{h} \sin^2(\pi m h),
  \ \mathrm{Im} \ \mu_m = \frac{1}{h} \sin(2\pi m h),
  \ m = 0, \pm 1, \pm 2, ..., \frac{M-1}{2}.
\end{equation}

\begin{figure}[h]
  \begin{center}
    \includegraphics[height=7.cm,angle=-0]{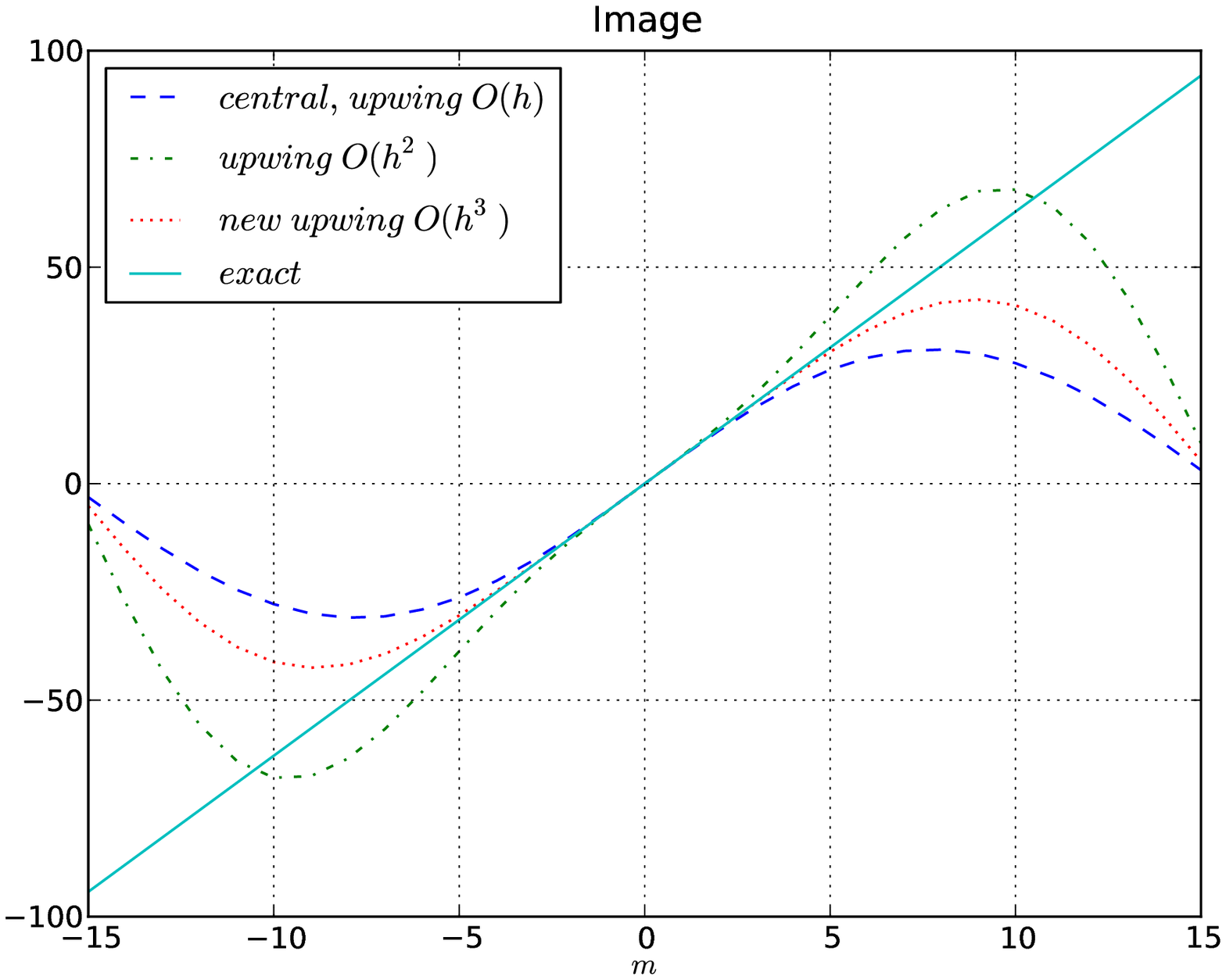} 
	\caption{}
	\label{f-1}
  \end{center}
\end{figure}

\begin{figure}[h]
  \begin{center}
    \includegraphics[height=7.cm,angle=-0]{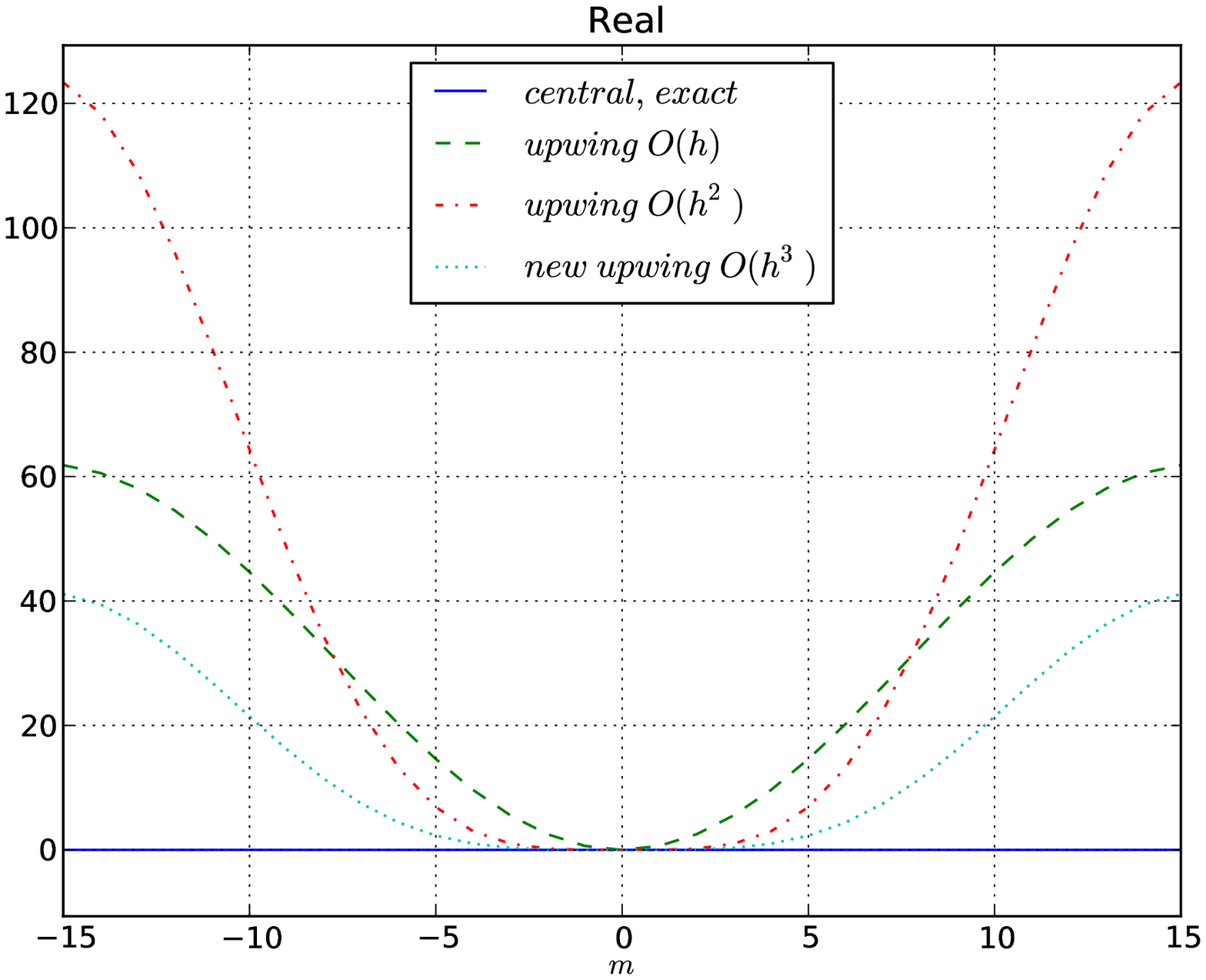} 
	\caption{}
	\label{f-2}
  \end{center}
\end{figure}

For central-difference approximations (\ref{15}) we obtain  
\begin{equation}\label{25}
  \mathrm{Re} \ \mu_m = 0,
  \quad \mathrm{Im} \ \mu_m = \frac{1}{h} \sin(2\pi m h),
  \quad m = 0, \pm 1, \pm 2, ..., \frac{M-1}{2}.
\end{equation}
Comparing (\ref{24}) and (\ref{25}) shows that the imaginary components of the 
spectrum of central-difference approximations 
and upwind difference approximations coincide.  
For upwind approximations we have positive real parts of the spectrum, which 
cause the dissipative properties of such approximations.  

Dissipative properties demonstrate also approximations(\ref{17}), (\ref{18}).
For the upwind second-order approximations we have  
\begin{equation}\label{26}
  \mathrm{Re} \ \mu_m = \frac{1}{h} (\cos(2\pi m h) - 1)^2,
  \quad \mathrm{Im} \ \mu_m = \frac{1}{h} \sin(2\pi m h)(2-\cos(2\pi m h))
\end{equation}
with the above-mentioned values of $ m $.  
For the third order approximations it is easy to obtain  
\begin{equation}\label{27}
  \mathrm{Re} \ \mu_m = \frac{1}{3 h} (\cos(2\pi m h) - 1)^2,
  \quad \mathrm{Im} \ \mu_m = \frac{1}{3 h} \sin(2\pi m h)(4-\cos(2\pi m h)) 
\end{equation}
respectively.

An illustration of the spectrum of grid convection  operator (\ref{24})--(\ref{27}) is shown in 
Fig.\ref{f-1},\ref{f-2} for  $M=31$. 
In particular, the main disadvantage of the scheme with directional differences of first order is associated 
with substantial dissipation of low harmonics whereas dissipative properties of scheme  (\ref{17}), (\ref{18}) 
are connected primarily with high harmonics.  
Approximation of  third-order for (\ref{18}) is relatively well reflects the spectral properties of the differential problem.
Its dissipative properties work only for high harmonics and is weak for the most important low harmonics of the difference solution.  

A similar analysis has been performed for the diffusion equation, where $\chi = 1$, $\mathcal{L} = \mathcal{D}$.
In this case, we investigated the non-negativity and self-adjointness of discrete diffusion operator $C$ and its spectral properties.  

\textit{1. Approximation of second order.}
The standard approximation at the three-point stencil leads us to  
\begin{equation}\label{28}
  D = - \partial_+ \partial_- .
\end{equation}

\textit{2. Approximation of fourth order.}
At the extended stencil we can use  
\begin{equation}\label{29}
  D = - \partial_+ \partial_- + 
  \frac{h^2}{12} \partial_+ \partial_- \partial_+ \partial_-.
\end{equation}
Using (\ref{28}), (\ref{29}), we have $D = D^* \geq 0$ in $H$.  

The spectrum of these operators $D$ is real with the same eigenfunctions that for $C$.  
For the eigenvalues we have  
\[
  \mu_m = \frac{4}{h^2} \sin^2 \frac{m \pi}{M},
  \quad m = 0, 1, ..., M-1
\]
if  (\ref{28}) and 
\[
  \mu_m = \frac{4}{h^2} \sin^2 \frac{m \pi}{M}
  \left (1 + \frac{1}{3} \sin^2 \frac{m \pi}{M} \right ),
  \quad m = 0, 1, ..., M-1
\]
if approximation (\ref{29}) are selected.
For the differential operator we have $\lambda_m = 4 \pi^2 m^2, \ m = 0, 1, ...$.
As expected, approximation (\ref{29}) gives us the better approximations for  the spectrum of the differential operator of diffusive transport.  

\section{SM properties of the approximation in time}

We'll use the two-level difference schemes for the approximate solution of  (\ref{10}), (\ref{11}).
Define a uniform grid in time with the time-step $\tau$ 
\[
  \overline{\omega}_\tau =
  \omega_\tau\cup \{T\} =
  \{t_n=n\tau,
  \quad n=0,1,...,N,
  \quad \tau N=T\} 
\]
and let $y_n = y(t_n), \ t_n = n \tau$.
For the exact solution of problem(\ref{10}), (\ref{11}), at transition from time level $t_n$  to  new time level $t_{n+1}$ we have  
\begin{equation}\label{30}
  y(x,t_{n+1}) =  e^{- \Lambda \tau} y(x,t_{n}) =
  \sum_{m=-(M-1)/2}^{(M-1)/2} (y(x,t_{n}), w_m) e^{-\mu_m \tau} w_m(x) .
\end{equation}

Two-level  difference scheme for problem  (\ref{10}), (\ref{11})  is written in the canonical operator-difference form 
\begin{equation}\label{31}
  B\frac{y_{n+1} -y_n}{\tau} +Ay_n = 0, \quad  n = 0,1, ... 
\end{equation}
with some operators $A$ and $B$.
In the Samarskii  theory of stability of operator-difference schemes \cite{SamTheor,SamVabMat,SamGul} 
stability conditions in the various norms are formulated in the form of operator inequalities for the $A, B$.  

Difference scheme (\ref{31})  is written as follows 
\begin{equation}\label{32}
y_{n+1} = S y_{n}, \qquad  n = 0,1, ... ,
\end{equation}
where  
\begin{equation}\label{33}
 S = E -\tau B^{-1} A
\end{equation}
is the  operator of transition from one time level to another level, which, in general, may depend on $n$. 

We restrict ourselves to the simplest difference approximation in time for problem (\ref{10}), (\ref{11}), 
which lead to the transition operator  
\begin{equation}\label{34}
 S = s(\tau \Lambda ),
\end{equation}
where $s(z)$ is a function of stability \cite{Hairer,Butcher}.
With constraints (\ref{34}) ($\tau A = (\tau A)(\tau \Lambda)$,
$B=B(\tau \Lambda)$)  the stability conditions in Hilbert spaces are easily verified on the basis of only  properties of function $s(z)$. 
Let $\Lambda \geq \delta E$, then  
\[
  \|s(\tau \Lambda ) \| \leq \max\limits_{\mathrm{Re} \, z \geq \delta \tau} |s(z) | ,
\]
and self-adjointness of operator  $\Lambda$ is not assumed.  

In the case of (\ref{34}) for the approximate solution at the new time level we have the representation 
\begin{equation}\label{35}
  y(x,t_{n+1}) =  \sum_{m=-(M-1)/2}^{(M-1)/2} (y(x,t_{n}), w_m) s(\mu_m \tau) w_m(x) .
\end{equation}
Quality of difference approximations in  time is estimated by comparing (\ref{35}) with representation (\ref{30}) for model problem (\ref{10}),(\ref{11}).
The comparison is performed at the level of behavior of individual harmonics and so we are talking about 
the SM properties for approximation in time.  

For convection problem  ($\chi = 1$, $\mathcal{L} = \mathcal{C}$) 
the spectrum is purely imaginary, and the solution is neutrally stable.  
After approximation in space using the above directional differences a typical situation is where the imaginary part of the spectrum is complemented by real part.  
When choosing approximations in time for the considered problems with skew-symmetric operators, we must monitor the behavior of the main imaginary part of the spectrum.  
This means that in problem  (\ref{10}), (\ref{11})
\[
  \Lambda = \Lambda_0 + \Lambda_1,
  \quad \Lambda_0 =\Lambda_0^* = \frac{1}{2} (\Lambda + \Lambda^*),
  \quad \Lambda_1 =- \Lambda_1^* = \frac{1}{2} (\Lambda - \Lambda^*)
\]
operator $\Lambda_1$  is main in the sense that $\Lambda_0 y \rightarrow 0$ as $h \rightarrow 0$ for sufficiently smooth  $y$. 
The supporting real part of the spectrum associated with operator $\Lambda_0$, is generated by the approximations in space and 
plays a minor role in these problems (it should not lead to instability of the difference solution).  

The difference scheme for convection problem  (\ref{10}), (\ref{11}),  
in which operator $\Lambda \geq 0$ and its antisymmetric part has the major role, 
is called as the SM stable if the difference scheme is stable and neutrally stable  at $\Lambda = - \Lambda^*$.

Two-level difference schemes of higher order accuracy for time-dependent linear problems we will construct on the basis of the Pad$\acute{e}$
approximations for the operator (matrix) exponent $e^{- \Lambda \tau}$.
For $e^{-z}$ we have  
\[
  e^{-z} = R_{lm}(z) + O(z^{l+m+1}),
  \quad R_{lm}(z) \equiv \frac{P_{lm}(z)}{Q_{lm}(z)},
\]
where $P_{lm}(z)$  and $Q_{lm}(z)$ are polynomials of degree $l$ and $m$, respectively:  \[
  P_{lm}(z) = \frac{l!}{(l+m)!} \sum_{k=0}^{l} 
  \frac{(l+m-k)!}{k! (l-k)!} (-z)^k,
\]
\[
  Q_{lm}(z) = \frac{m!}{(l+m)!} \sum_{k=0}^{m} 
  \frac{(l+m-k)!}{k! (m-k)!} z^k .
\]
For equation (\ref{10}) the application of Pad$\acute{e}$ approximations corresponds to the two-level scheme  
\begin{equation}\label{36}
  Q_{lm}(\tau \Lambda) \frac{y_{n+1} -y_n}{\tau} + 
  \frac{1}{\tau }(Q_{lm}(\tau \Lambda) - P_{lm}(\tau \Lambda)) y_{n} = 0, 
  \quad  n = 0,1, ... .
\end{equation}
In canonical form (\ref{31}) difference scheme (\ref{36}) corresponds to the choice  
\begin{equation}\label{37}
  A = \frac{1}{\tau }(Q_{lm}(\tau \Lambda) - P_{lm}(\tau \Lambda)),
  \quad B = Q_{lm}(\tau \Lambda) .
\end{equation}

Difference schemes for problem (\ref{10}), 
(\ref{11}) with  $\Lambda \geq 0$
on the basis of Pad$\acute{e}$ approximations are stable (absolutely stable) (estimate $\|y_{n+1}\| \leq \|y_{n}\|$ holds) 
at $l \leq m$ \cite{Hairer,Butcher}. 
It is only necessary  to highlight among such schemes the SM-stable schemes.  

In the simplest case $m = 1$  we have  
\[
  R_{01}(z) = \frac{1}{1 + z} = e^{-z} + O(z^{2}),
\]
\[
  R_{11}(z) = \frac{1-\frac12 z }{1 + \frac12 z} = 
  e^{-z} + O(z^{3}) .
\]
Approximations  $R_{01}(z)$ corresponds to application for the approximate solution of problem (\ref{10}), (\ref{11}) 
the purely implicit scheme 
\begin{equation}\label{38}
  \frac{y_{n+1} -y_n}{\tau} + \Lambda y_{n+1} = 0, \quad  n = 0,1, ... .
\end{equation}
The application of the symmetric scheme (Crank-Nicholson)  
\begin{equation}\label{39}
  \frac{y_{n+1} -y_n}{\tau} + \Lambda \frac{y_{n+1}+y_{n}}{2} = 0, 
  \quad  n = 0,1, ... .
\end{equation}
corresponds to the choice of approximation  $R_{11}(z)$.

The condition of neutral stability of $\|y_{n+1}\| = \|y_{n}\|$ for two-level scheme (\ref{32})  will be satisfied at $\|S\| = 1$.
Taking into account  (\ref{34}) and for $\Lambda = - \Lambda^*$ this corresponds to the case   
\begin{equation}\label{40}
  |s(z)| = |R_{lm}(z)| = 1, \quad \mathrm{Re} \ z = 0.
\end{equation}

For purely implicit scheme (\ref{38}) we have  
\[
  |R_{01}(z)| = \frac{1}{\sqrt{1 + y^2}},
  \quad z = i y.
\]
Thus, the condition of neutral stability is not satisfied --- the purely implicit is not SM stable for problems 
with the main skew-symmetric operator.    
While for  symmetric scheme (\ref{39}) we obtain  
\[
  |R_{11}(z)| = 1,
  \quad z = i y.
\]
Thus, this scheme is SM stable for the investigated class of problems.

We can make similar conclusions for schemes with Pad$\acute{e}$ approximations at $m > 1$. 
Only a scheme that is based on approximation $R_{mm}$ is SM stable for problems with the main skew-symmetric  operator.  
At using Pad$\acute{e}$ approximations with $l <m$ the scheme demonstrates dissipative properties due to the approximation in time.  
Only at $l = m$ the corresponding scheme is neutrally stable.  

A similar analysis is carried out (see\cite{VabSM}) for the diffusion equation.  
In  problem (\ref{10}), (\ref{11}) with $\Lambda = \Lambda^* \geq 0$ amplitudes of harmonics with higher numbers damp  
more quickly in compare with amplitudes of harmonics with lower numbers (spectral monotonicity) and damp  to zero as $t \rightarrow \infty$ 
(asymptotic stability).  Such a behavior of the approximate solutions we associate with the SM properties of approximation in time 
for the solution of  problems with self-adjoint operators.  

We assume that the difference scheme for problem (\ref{10}), (\ref{11}) with $\Lambda = \Lambda^* \geq 0$ is  SM stable 
if it is spectrally monotonic and asymptotically stable.  

Difference schemes based on the Pad$\acute{e}$ approximation $R_{lm}$ are SM stable at $ l = 0 $.
Purely implicit scheme (\ref{38}) belongs to this class of schemes, 
whereas the symmetric scheme is conditionally SM stable.  

The main conclusion of the study is that the approximate solution of problems with skew-symmetric operators 
we must use such approximations in time, which are based on the  Pad$\acute{e}$ approximations $R_{mm}(z)$. 
For problems with self-adjoint operators it is more preferred to use Pad$\acute{e}$ approximation $R_{0m}(z)$.
For problems with general non-selfadjoint  operators, approximation in time can be constructed via decomposition into 
the self-adjoint and skew-symmetric  components and further constructing different approximations for them, 
based on special splitting schemes.


\begin{thebibliography}{1}

\bibitem{Hundsdorfer}
\newblock {\em Hundsdorfer W., Verwer J.}
Numerical Solution of Time-Dependent
  Advection-diffusion-reaction Equations. Berlin: Springer, 2003.

\bibitem{Gustafsson}
\newblock {\em Gustafsson B.}
High Order Difference Methods for Time Dependent PDE. Berlin: Springer, 2008.

\bibitem{Ascher}
\newblock {\em Ascher U.~M.}
Numerical Methods for Evolutionary Differential Equations.
Philadelphia, PA: Society for Industrial and Applied Mathematics (SIAM), 2008.

\bibitem{LeVeque}
\newblock {\em LeVeque R.~J.}
Finite Difference Methods for Ordinary and Partial Differential
Equations. Steady-state and Time-dependent problems.
Philadelphia, PA: Society for Industrial and Applied Mathematics
(SIAM), 2007.

\bibitem{Rakit}
\newblock {\em Rakitskii Yu. V., Ustinov S.~M.,  
Chernorutskii I. G.}
Numerical Methods for Solving Stiff Systems.
Moscow: Nauka, 1979, in Russian.

\bibitem{Hairer}
\newblock {\em Hairer E., Wanner G.}
Solving Ordinary Differential Equations. II: Stiff and
Differential-Algebraic Problems. Berlin: Springer, 1996.

\bibitem{Butcher}
\newblock {\em Butcher J.~C.}
Numerical Methods for Ordinary Differential Equations. 
Hoboken, NJ: Wiley, 2008.

\bibitem{Verwer}
\newblock {\em Dekker K., Verwer J.}
Stability of Runge-Kutta Methods for Stiff Nonlinear
Differential Equations. Amsterdam - New York - Oxford: North-Holland, 1984.

\bibitem{Gear}
\newblock {\em Gear C.~W.}
Numerical Initial Value Problems in Ordinary Differential
Equations. Englewood Cliffs, NJ: Prentice-Hall, 1971.

\bibitem{SamTheor}
\newblock {\em Samarskii, A.~A.}
The Theory of Difference Schemes. New York:
Marcel Dekker Inc.,  2001.
b
\bibitem{SamVabMat}
\newblock {\em Samarskii A.~A., Matus P.~P., Vabishchevich P.~N.}
Difference Schemes with Operator Factors. Dordrecht Hardbound:
Kluwer Academic Publishers, 2002.

\bibitem{SamGul}
\newblock {\em Samarskii A.A., Gulin A.~V.}
Stability of Difference Schemes.
Moscow: Nauka, 1973, in Russian.

\bibitem{SamVab}
\newblock {\em Samarskii A.~A.,  Vabishchevich P.~N.}
Computational Heat Transfer. Vol. 1. Mathematical Modelling. Chichester:
Wiley, 1995.

\bibitem{VabSM}
\newblock {\em Vabishchevich P.~N.}
Two-Level Finite Difference Scheme of Improved Accuracy
Order for Time-Dependent Problems of Mathematical Physics.
\newblock {\em Computational Mathematics and Mathematical Physics}, 2010, 
Vol. 50, No. 1, pp. 112-123.

\bibitem{SamVabConv}
\newblock {\em Samarskii A.~A.,  Vabishchevich P.~N.}
Methods for Convection-Diffusion Problems. 
Moscow: URSS, 2004, in Russian.

\bibitem{morton}
\newblock {\em Morton K.~W.}
Numerical Solution of Convection-Diffusion Problems.
London: Chapman \& Hall, 1996.

\bibitem{Hirsch}
\newblock {\em Hirsch C.}
Numerical Computation of Internal and External Flows.
Fundamentals of Computational Fluid Dynamics.
Amsterdam: Butterworth-Heinemann, 2007.

\end{thebibliography}
\end{document}